\documentclass[11pt]{article}
\UseRawInputEncoding
\usepackage{amsfonts}
\usepackage{latexsym,amsmath}
\usepackage{amssymb,array}
\usepackage[sort&compress,round,comma]{natbib}
\usepackage{setspace}
\usepackage{placeins}
\makeatletter \oddsidemargin 0in \evensidemargin 0in \textwidth 16cm
\RequirePackage[dvips]{graphicx}  \textheight 18cm
\newtheorem{Definition}{Definition}
\date{}

\begin{document}
\title{A multi-objective reliability-redundancy allocation problem with active redundancy and interval type-2 fuzzy parameters}
\author{Pradip Kundu\footnote{Corresponding author, e-mail: kundu.maths@gmail.com}\\\textit{\small{Decision Science and Operations Management, School of Management,}}\\
\textit{\small{Birla Global University, Bhubaneswar, Odisha-751003,
India}}}
\maketitle
\begin{abstract}
This paper considers a multi-objective reliability-redundancy
allocation problem (MORRAP) of a series-parallel system, where
system reliability and system cost are to be optimized
simultaneously subject to limits on weight, volume, and redundancy
level. Precise computation of component reliability is very
difficult as the estimation of a single number for the probabilities
and performance levels are not always possible, because it is
affected by many factors such as inaccuracy and insufficiency of
data, manufacturing process, environment in which the system is
running, evaluation done by multiple experts, etc. To cope with
impreciseness, we model component reliabilities as interval type-2
fuzzy numbers (IT2 FNs), which is more suitable to represent
uncertainties than usual or type-1 fuzzy numbers. To solve the
problem with interval type-2 fuzzy parameters, we first apply
various type-reduction and defuzzification techniques, and obtain
corresponding defuzzified values. As maximization of system
reliability and minimization of system cost are conflicting to each
other, so to obtain compromise solution of the MORRAP with
defuzzified parameters, we apply five different multi-objective
optimization methods, and then corresponding solutions are analyzed.
The problem is illustrated numerically for a real-world MORRAP on
pharmaceutical plant, and solutions are obtained by standard
optimization solver LINGO, which is based on gradient-based
optimization - Generalized Reduced Gradient (GRG) technique.
\end{abstract}
\section{Introduction}
An industrial or a mechanical system such as aircraft, nuclear
plants, lighting system, material handling systems, pharmaceutical
plant, civil engineering systems, and so on is composed of numerous
complex components. The reliability, i.e., the probability that a
system performs satisfactorily over a certain period of time depends
on each of its constituent components and the system design. The
study of reliability optimization relates to enhance the reliability
of a system so that the system can be operational satisfactorily for
the maximum possible time. Reliability of a system can be improved
by using high reliable components and adding redundant (standby)
components. However, this may increase the system cost. Further
designing of system structure also depends on various
resource/engineering constraints related to cost, volume, weight,
and energy consumption, etc. Reliability-redundancy allocation
problem (RRAP) (\citealp{kuo}) is the problem of maximizing system
reliability through redundancy (\citealp{case}) and component
reliability choices subject to the applicable resource constraints.
However, in addition to maximization of system reliability, if the
system cost or weight has to be minimized simultaneously, then the
problem becomes the multi-objective reliability-redundancy
allocation problem (MORRAP) (\citealp{arda}; \citealp{cao};
\citealp{garg}; \citealp{huan}; \citealp{khal}; \citealp{muhu};
\citealp{rao}; \citealp{safa}). The main goal of such MORRAP is to
determine the optimal component reliabilities and number of
redundant components in each of the subsystems to maximize the
system reliability and minimize the system cost simultaneously
subject to several resource constraints.

Practically, exact computation of reliability is challenging and is
associated with uncertainties due to various reasons. Notably, the
estimation of a single number for the probabilities and performance
levels is very difficult (\citealp{cheng}). Some reasons come from
inaccuracy and insufficiency of data, data collection from multiple
sources or evaluation done by multiple experts, etc. Some other key
sources of uncertainties are uncertainty factors in manufacturing
process (like quality assurance controls, work management and
execution, maintenance activities) and environmental factors (the
reliability of components depends on the factors like temperature,
voltage, humidity of the environment in which the associated system
is running). So, it is not always possible to precisely determine
the reliability of a component. Many researchers have investigated
RRAP with various uncertain parameters such as interval-valued
(\citealp{saho}; \citealp{garg1}; \citealp{roy}; \citealp{xu};
\citealp{zhan}), and fuzzy parameters (\citealp{maha};
\citealp{yao}; \citealp{cheng}; \citealp{garg2}; \citealp{kuma};
\citealp{srir}; \citealp{muhu}). Most of the RRAP with fuzzy
parameters considered \textit{type-1 fuzzy numbers} (T1 FNs) for
modeling the uncertainties (\citealp{chen}; \citealp{alie};
\citealp{yao}; \citealp{jamk}).

\textit{Type-2 fuzzy sets} (T2 FSs) (for details see
\citealp{mend5}) are proved to be more suitable in many instances to
represent uncertainties than ordinary or \textit{type-1 fuzzy sets}.
A \textit{Type-2 fuzzy set} is a generalization of \textit{type-1
fuzzy set}, and has an extra degree of freedom to represent
uncertainties because of its secondary membership function.
\textit{Interval type-2 fuzzy set} (IT2 FS) (\citealp{mend6}), a
special case of general T2 FS has been successfully used to model
uncertain parameters in many instances like data collected from
multiple sources, opinion taken from several experts, information
given by approximate intervals or linguistic terms, etc.
Specifically, it is showed by many researchers that to deal with
linguistic uncertainties (\citealp{mend2,mend3,mend4};
\citealp{liu3}; \citealp{mill}), approximate intervals (like two
endpoints of the intervals are not certain) (\citealp{liu3}) or
several membership functions (\citealp{pago}), \textit{interval
type-2 fuzzy set} is an appropriate tool. Not only that,
\cite{mend2,mend3,mend4} explained and showed that modeling
linguistic information using \textit{type-1 fuzzy set} is not
scientific, instead one should use \textit{type-2 fuzzy set},
specifically \textit{interval type-2 fuzzy set}. As we mention
earlier, there are several research work has been done on RRAP with
T1 FNs. However, there are very few research works on RRAP with
\textit{type-2 fuzzy numbers} available in the literature
(\citealp{muhu}). The significant contributions of the present
investigation are as follows:
\begin{itemize}
\item We formulate a MORRAP of a series-parallel system with the approximate reliability of each component of a subsystem represented as \textit{interval type-2 fuzzy numbers} (IT2 FNs). Most of the previous research work has been investigated RRAP with interval numbers or T1
FNs.
\item We not only explain but also illustrate numerically that modeling uncertain parameters (reliabilities) using IT2 FNs leads to the better performance than that of using T1 FNs, i.e. our investigation suggest that we can model system with higher system reliability and less system cost.
\item We apply various type-reduction and defuzzification techniques to obtain corresponding defuzzified values of IT2 FNs, and comparative study has been presented.
\item To deal with conflicting objectives we apply five different multi-objective optimization techniques to obtain solution of the problem. As a result a decision maker can choose appropriate result according to his/her preference or as situation demand.
\end{itemize}

In our considered MORRAP there are two conflicting objectives,
namely, maximization of system reliability and minimization of
system cost. Construction of IT2 FNs to represent imprecise
component reliabilities has been done by using a modified algorithm
which was initially proposed by \cite{muhu}. To solve MORRAP with
interval type-2 fuzzy parameters, we first apply various
type-reduction and defuzzification techniques to obtain
corresponding defuzzified values. To deal with two conflicting
objectives we then apply five different multi-objective optimization
methods, and obtain compromise solution of the problem. The problem
is also solved by modeling component reliabilities as T1 FNs, and
the obtained result is compared with the result for the same problem
with IT2 FNs. The rest of the paper is organized as follows. Section
2 provides brief introduction of \textit{type-2 fuzzy set}. The
detail of the problem (MORRAP) formulation is presented in Section
3. Section 4 discusses some type-reduction and defuzzification
techniques in brief. Section 5 presents some multi-objective
optimization techniques in detail. In Section 6, the problem and
methods are illustrated numerically for a real-world MORRAP on
pharmaceutical plant. Finally, Section 7 concludes the paper.

\section{Preliminaries}
\subsection{Type-2 fuzzy set}
\textit{Type-2 fuzzy set} (T2 FS) is an extension of usual or
\textit{type-1 fuzzy set} (T1 FS). It is a fuzzy set with fuzzy
membership function, i.e., membership grade of each element in the
set is no longer a precise (crisp) value but a fuzzy set.
\begin{Definition} A type-2 fuzzy set $\tilde{\tilde{A}}$ (\citealp{mend5}) in
a space of points (objects) $X$ is characterized by a type-2
membership function $\mu_{\tilde{\tilde{A}}}:X\times
J_{x}\rightarrow [0,1]$, and is defined as
 $$\tilde{\tilde{A}}=\{((x,u),\mu_{\tilde{\tilde{A}}}(x,u)):~\forall x\in X,~\forall u\in J_{x}\subseteq [0,1]\},$$
where $J_{x}\subseteq [0,1]$ is the primary membership of $x\in X$,
and $0\leq \mu_{\tilde{\tilde{A}}}(x,u)\leq 1$ for all $x\in X$,
$u\in J_{x}$. $\tilde{\tilde{A}}$ is also expressed as
$$\tilde{\tilde{A}}=\int_{x\in X} \int_{u\in J_{x}}
\mu_{\tilde{\tilde{A}}}(x,u)/(x,u)~,~J_{x}\subseteq [0,1],$$ where
$\int\int$ denotes union over all admissible $x$ and $u$. For
discrete universes of discourse, $\int$ is replaced by $\sum$.\\ For
particular $x=x'\in X$, $\mu_{\tilde{\tilde{A}}}(x',u)$ $\forall
u\in J_{x'}$, is called secondary membership of $x'$. The amplitude
of a secondary membership function is called a secondary membership
grade. Thus $\mu_{\tilde{\tilde{A}}}(x',u')$, $u'\in J_{x'}$ is
secondary membership grade of $(x',u')$ which represents the grade
of membership that the point $x'$ has the primary membership
$u'$.\end{Definition}
\begin{Definition} An interval type-2 fuzzy set (IT2 FS) (\citealp{mend6}) is a special case of
T2 FS where all the secondary membership grades are 1, i.e.,
$\mu_{\tilde{\tilde{A}}}(x,u)=1$ for all $(x,u)\in X\times J_x$. An
IT2 FS $\tilde{\tilde{A}}$ can be written as
$$\tilde{\tilde{A}}=\int_{x\in X} \int_{u\in J_{x}}
1/(x,u)=\int_{x\in X} \left. \left[\int_{u\in J_{x}}
1/u\right]\middle/ x \right.,~J_{x}\subseteq [0,1].$$ As the
secondary membership grades are 1, an IT2 FS can be characterized by
the footprint of uncertainty (FOU) which is the union of all primary
memberships $J_{x}$ in a bounded region, so that it is defined as $$
FOU(\tilde{\tilde{A}}) = \bigcup_{x\in X}~J_{x}.$$ \end{Definition}
The FOU (see Fig. 1) is bounded by an upper membership function
(UMF) $\bar{\mu}_{\tilde{\tilde{A}}}(\cdot)$ and a lower membership
function (LMF) $\underline{\mu}_{\tilde{\tilde{A}}}(\cdot)$, both of
which are the membership functions of T1 FSs, and
$J_{x}=[\underline{\mu}_{\tilde{\tilde{A}}}(x),\bar{\mu}_{\tilde{\tilde{A}}}(x)]$,
$\forall~x\in X$. In this view, the IT2 FS can be represented by
$(\tilde{A^U},\tilde{A^L})$, where $\tilde{A^U}$ and $\tilde{A^L}$
are T1 FSs. The support of IT2 FS $\tilde{\tilde{A}}$ can written as
$supp(\tilde{\tilde{A}})=\{x\in X:
\bar{\mu}_{\tilde{\tilde{A}}}(x)>0\}$.
\begin{Definition} Interval type-2 fuzzy number (IT2 FN): An interval type-2 fuzzy
number (IT2 FN) (\citealp{hesa}) is an IT2 FS on set of real numbers
$\mathbb{R}$, whose upper and lower membership functions are
membership functions of T1 FNs.
\end{Definition} For example, Fig. 1 represents a triangular IT2 FN
$\tilde{\tilde{A}}=(\tilde{A^U},\tilde{A^L})=((2,4,6),(3,4,5))$,
where $\tilde{A^U}$ and $\tilde{A^L}$ are triangular fuzzy numbers
having following membership functions:
$$\mu_{\tilde{A^U}}(x)=\bar{\mu}_{\tilde{\tilde{A}}}(x)=\left\{
                        \begin{array}{ll}
                          \frac{x-2}{2}, & \hbox{if $2\leq x\leq 4$;} \\
                          1, & \hbox{if $x=4$;} \\
                          \frac{6-x}{2}, & \hbox{if $4\leq x\leq 6$;} \\
                          0, & \hbox{otherwise.}
                        \end{array}
                      \right.~\mu_{\tilde{A^L}}(x)=\underline{\mu}_{\tilde{\tilde{A}}}(x)=\left\{
                        \begin{array}{ll}
                          x-3, & \hbox{if $3\leq x\leq 4$;} \\
                          1, & \hbox{if $x=4$;} \\
                          5-x, & \hbox{if $4\leq x\leq 5$;} \\
                          0, & \hbox{otherwise.}
                        \end{array}
                      \right.$$
The arithmetic operations between two triangular IT2 FNs
$\tilde{\tilde{A_1}}=(\tilde{A_{1}^U},\tilde{A_{1}^L})=((a_{11}^U,a_{12}^U,a_{13}^U),(a_{11}^L,a_{12}^L,a_{13}^L))$
and
$\tilde{\tilde{A_2}}=(\tilde{A_{2}^U},\tilde{A_{2}^L})=((a_{21}^U,a_{22}^U,a_{23}^U),(a_{21}^L,a_{22}^L,a_{23}^L))$
are defined as follows:\\
Addition operation: $\tilde{A_{1}}\oplus
\tilde{A_{2}}=(\tilde{A_{1}^U},\tilde{A_{1}^L})\oplus
(\tilde{A_{2}^U},\tilde{A_{2}^L})$\\
$=((a_{11}^U+a_{21}^U,a_{12}^U+a_{22}^U,a_{13}^U+a_{23}^U),
(a_{11}^L+a_{21}^L,a_{12}^L+a_{22}^L,a_{13}^L+a_{23}^L)),$\\
Multiplication operation: $\tilde{A_{1}}\otimes
\tilde{A_{2}}=(\tilde{A_{1}^U},\tilde{A_{1}^L})\otimes
(\tilde{A_{2}^U},\tilde{A_{2}^L})$\\
$=((a_{11}^U\times a_{21}^U,a_{12}^U\times a_{22}^U,a_{13}^U\times
a_{23}^U), (a_{11}^L\times a_{21}^L,a_{12}^L\times
a_{22}^L,a_{13}^L\times
a_{23}^L)).$\\
The arithmetic operations between triangular IT2 FN
$\tilde{\tilde{A_1}}$ and a real number $r$ are defined as follows:\\
$r\tilde{\tilde{A_1}}=((r\times a_{11}^U,r\times a_{12}^U,r\times
a_{13}^U),(r\times a_{11}^L,r\times
a_{12}^L,r\times a_{13}^L)),$\\
$\frac{\tilde{\tilde{A_1}}}{r}=((\frac{1}{r}\times
a_{11}^U,\frac{1}{r}\times a_{12}^U,\frac{1}{r}\times
a_{13}^U),(\frac{1}{r}\times a_{11}^L,\frac{1}{r}\times
a_{12}^L,\frac{1}{r}\times a_{13}^L))$, where $r> 0$.
\begin{figure}\label{figit2fn}\begin{center}
  \includegraphics[width=10cm]{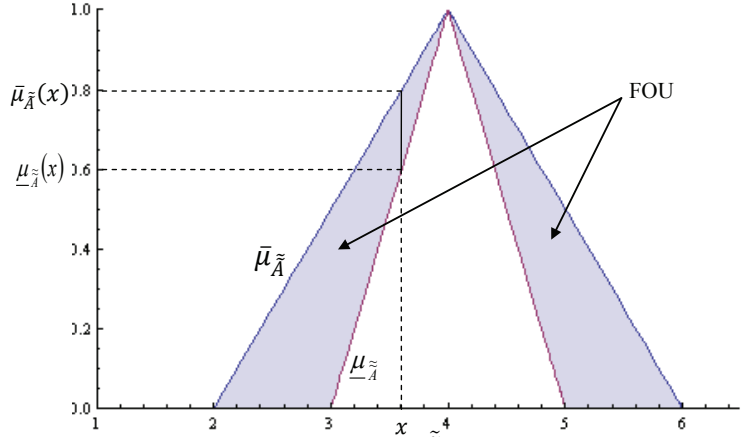}\\\vspace{-0.6cm}
  \caption{Triangular IT2 FN $\tilde{\tilde{A}}=((2,4,6),(3,4,5))$.}\end{center}
\end{figure}
\section{A Multi-objective reliability-redundancy allocation
problem (MORRAP)} Generally, complex systems are composed of several
subsystems (stages), each having more than one component. In
reliability context, system designing mainly concern with
improvement of overall system reliability, which may be subject to
various resource/engineering constraints associated with system
cost, weight, volume, and energy consumption. This may be done (i)
by incorporating more reliable components (units) and/or (ii)
incorporating more redundant components. In case of the second
approach, optimal redundancy is mainly taken into consideration for
the economical design of systems. Again the reliability optimization
concerned with redundancy allocation is generally classified into
two categories: (i) maximization of system reliability subject to
various resource constraints; and (ii) minimization of system cost
subject to the condition that the associated system reliability is
required to satisfy a desired level. However, if maximization of
system reliability and minimization of the system cost have to be
done simultaneously, then the problem becomes the multi-objective
reliability-redundancy allocation problem (MORRAP). So, the main
goal of MORRAP is to determine the optimal component reliabilities
and number of redundant components in each of the subsystems to
maximize the system reliability and minimize the system cost
simultaneously subject to several resource constraints.

Here, we have considered a MORRAP for a series-parallel system
configuration (\citealp{huan}; \citealp{garg}). A series-parallel
system usually has $m$ (say) independent subsystems arranged in
series, and in each subsystem, there are $n_i$ (say) $(i=1,2,...,m)$
components, which are arranged in parallel. A reliability block
diagram (RBD) of this series-parallel system is depicted in Fig. 2,
where small rectangular blocks represent the components in each of
the $m$ subsystems. The reliability block diagram provides a
graphical representation of the system that can be used to analyze
the relationship between component states and the success or failure
of a specified system. As seen from Fig. 2, in each subsystem the
components are arranged in parallel, so each of the subsystems can
work if at least one of its components works. Again as these
subsystems are arranged in series, the whole system can work if all
the subsystems work. Obviously, reliability of the series-parallel
system is the product of all the associated subsystem reliabilities.
For the considered MORRAP, the objective functions are maximization
of system reliability and minimization of system cost, subject to
limits on weight, volume, and redundancy level. Also, the problem
considers the active redundancy strategy (i.e., all the components
in each subsystem are active and arranged in parallel).
\begin{figure}\begin{center}
  \includegraphics[width=10cm]{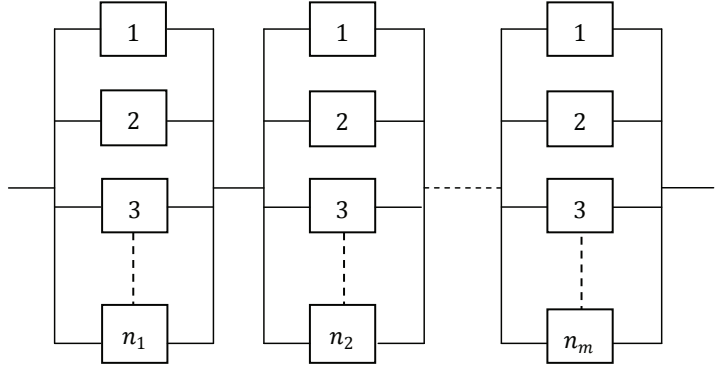}\\\vspace{-0.4cm}
  \caption{Reliability block diagram of series-parallel system.}\end{center}
\end{figure}

For the mathematical formulation of the problem we use the following
notations:
\begin{tabbing}
aaaaaaaa \= ababababababab \kill
  $m$ \>  Number of subsystems\\
  $n_i$ \> Number of components in subsystem $i$ ($i=1,2,...,m)$\\
  $r_i$\> Reliability of each of the components in $i$th subsystem\\
  $R_s$ \> System reliability\\
  $C_s$\> System cost \\
  $c(r_i)$\> Cost of the each component with reliability $r_i$ at subsystem
  $i$\\
  $\alpha_{i}$ \> Constants representing the physical characteristic
  (scaling factor) of the cost-reliability\\\> curve of each component with reliability $r_i$ at subsystem $i$\\
  $\beta_{i}$ \> Constants representing the physical characteristic
  (shaping factor) of the cost-reliability\\\> curve of each component with reliability $r_i$ at subsystem $i$\\
  $T$ \> Operating time during which the component must not fail \\
  $w_{i}$ \> Weight of each component of the subsystem $i$\\
  $W_{s}$ \> System weight\\
  $W$ \> Upper limit on the weight of the system\\
  $v_i$ \> Volume of each component of the subsystem $i$\\
  $V_{s}$ \> System volume\\
  $V$ \> Upper limit on the volume of the system
 \end{tabbing}
It is assumed that all the components for individual subsystem are
identical, all redundancies are active, failures of individual
components are independent, and each component can only be in one of
two states, i.e., either working or failed. Now the mathematical
formulation of the MORRAP is as follows:
\begin{eqnarray}\label{eqmorp1}\text{Max}~R_{s}&=&\prod_{i=1}^m~\left[1-(1-r_i)^{n_i}\right],\\
\label{eqmorp2}Min~C_s&=&\sum_{i=1}^m
c(r_i)(n_i+\exp(n_i/4)),\end{eqnarray}
\begin{eqnarray}\label{eqmorp3} \text{subject to}
&&V_{s}=\sum_{i=1}^m v_{i} n_{i}^2~\leq V,\\
\label{eqmorp4}&&W_{s}=\sum_{i=1}^m w_{i}(n_i\cdot \exp(n_i/4))\leq
W,
\end{eqnarray}
\begin{equation}\label{eqmorp5}r_{i,min}\leq r_i \leq r_{i,max}, 1\leq n_i \leq n_{i,max}, r_i\in(0,1), i=1,2,...,m.\end{equation}
For the presented model, cost of the each component is an increasing
function of the component reliability or conversely a decreasing
function of the failure rate (hazard rate) of the component, where
failure times of components follow exponential distribution. So the
reliability of each of the component in subsystem $i$,
\begin{equation}\label{ri}r_i=e^{-\lambda_i T},~\lambda_i>0,\end{equation} and
consequently the hazard rate is $\lambda_i$, where $T$ be the
operating time during which the component will not fail. As cost of
the each component in the $i$-th subsystem, $c(r_i)$, is a
decreasing function of the hazard rate, it is represented as
\begin{equation}\label{cri}c(r_i)=\alpha_{i} \lambda_i^{-\beta_{i}},\end{equation} where
$\alpha_{i}$ and $\beta_{i}$ are constants. Equations (\ref{ri}) and
(\ref{cri}) together gives
$c(r_i)=\alpha_{i}\left(\frac{-T}{ln(r_i)}\right)^{\beta_{i}}$. Now,
each subsystem is comprised of some components connected in
parallel. The factor $\exp(n_i/4)$ is incorporated due to the
interconnecting hardware between the parallel components
(\citealp{rao}; \citealp{pras}; \citealp{wang}; \citealp{arda}).
Total volume of the system $(V_{s})$ which consists of the volume of
the each component as well as space between the components and space
between the subsystems, is represented in equation (\ref{eqmorp3}).
Here $n_{i,max}$ represents the maximum number of components allowed
in subsystem $i$ arranged in parallel, and $r_{i,min}$ and
$r_{i,max}$ respectively the minimum and maximum reliability limits
of each component in subsystem $i$.
\subsection{MORRAP with interval type-2 fuzzy parameters}\label{secMOIT2}
As discussed in the introduction, the component reliability in a
system cannot be always precisely measured as crisp values, but may
be determined as approximate values like ``about 0.6" or approximate
intervals with imprecise end points. Some of the reasons are
inaccuracy and insufficiency of data, manufacturing uncertainty,
environmental issues (like temperature, humidity of the environment
in which the system is running), evaluation done by multiple experts
or data collected from multiple sources, etc. So to cope with the
ambiguity/approximation we associate a degree of membership to each
value of reliability. Here the approximate reliability of each
component of a subsystem is represented by IT2 fuzzy number and is
denoted by $\tilde{\tilde{r}}_i$, $i=1,2,...,m$. The assumption of
IT2 FN to represent uncertainty is vary reasonable when value of a
parameter is given by approximate interval (like the two end-points
of the interval are not exact), linguistic terms, etc. Now the above
MORRAP (\ref{eqmorp1})-(\ref{eqmorp5}) becomes
\begin{eqnarray}\label{morpt21}&&\text{Max}~\tilde{\tilde{R}}_{s}=\prod_{i=1}^m~\left[1-(1-\tilde{\tilde{r}}_i)^{n_i}\right],\\
\label{morpt22}&&Min~\tilde{\tilde{C}}_s=\sum_{i=1}^m
\alpha_{i}\left(\frac{-T}{ln(\tilde{\tilde{r}}_i)}\right)^{\beta_{i}}(n_i+\exp(n_i/4))\\\label{morpt23}
&&\text{subject to }~(\ref{eqmorp3})-(\ref{eqmorp5}).\end{eqnarray}
To solve this problem we use different type-reduction and
corresponding defuzzification strategies to convert the problem with
IT2 fuzzy parameters to the problem with defuzzified parameters.
Then we use various multi-objective techniques to solve the
deterministic bi-objective problem. To construct interval type-2
fuzzy membership function for the reliability $\tilde{\tilde{r}}_i$
having support $supp(\tilde{\tilde{r}}_i)\subseteq [a,b]\subset
[0,1]$ we use the following algorithm. To construct this algorithm
we modified the Algorithm-1 of \cite{muhu} to ensure that the
support of $\tilde{\tilde{r}}_i$ must lie within $(0,1)$.

\textbf{Algorithm: Generation of T1 FN $(\tilde{r}_i)$ and IT2 FN
$(\tilde{\tilde{r}}_i)$}
\begin{description}
  \item[Step 1:] Take $r_i\in [a,b]$.
  \item[Step 2:] Find the values of $r_{i}^l$ and $r_{i}^u$ as follows:\\
$r_{i}^l=a+(r_{i}-a)*rand$\\
$r_{i}^u=b-(b-r_{i})*rand$
\item[Step 3:] Construct the T1 FN
$\tilde{r}_i=(r_{i}^l,r_i,r_{i}^u)$.~// \textit{This step should be
skipped for generation of IT2 FN.}
  \item[Step 4:] Find the values of $r_{i1}^l$ and $r_{i3}^l$ as follows:\\
$r_{i1}^l=r_{i}^l+(r_{i}-r_{i}^l)*rand$\\
$r_{i3}^l=r_{i}^u-(r_{i}^u-r_{i})*rand$
\item[Step 5:] Find the values of $r_{i1}^u$ and $r_{i3}^u$ as follows:\\
$r_{i1}^u=r_{i}^l-(r_{i}^l-a)*rand$\\
$r_{i3}^u=r_{i}^u+(b-r_{i}^u)*rand$
\item[Step 6:] Construct the IT2 FN
$\tilde{\tilde{r}}_i=((r_{i1}^u,r_{i},r_{i3}^u),(r_{i1}^l,r_{i},r_{i3}^l))$.
\end{description}

In the next section we have briefly introduced different
type-reduction and defuzzification strategies of interval type-2
fuzzy set.
\section{Type-reduction and defuzzification
strategies}\label{sectype} Here, we discuss some type-reduction and
defuzzification strategies those are investigated in this study to
obtain corresponding type-reduced set and defuzzified values of
interval type-2 fuzzy parameters. These methods are given in detail
in the corresponding references. However, we present the methods
briefly to provide a ready reference to readers.
\subsection{Karnik-Mendel (KM) algorithm}
\cite{karn} introduced the concept of centroid of T2 FS by which it
can be reduced to a T1 FS (\citealp{liuf}). The computational
procedure to find centroid of an IT2 FS $\tilde{\tilde{A}}$ starts
with discretization (if the domain is not discrete) of the
continuous domain into finite number of points $x_i$, $i=1,2,...,N$
which are sorted in an ascending order. Then the centroid of the IT2
FS is given by $[y_l,y_r]$ and corresponding defuzzified value is
$C(\tilde{\tilde{A}})=(y_l+y_r)/2$, where
$$y_l=\frac{\sum_{i=1}^L x_i \bar{\mu}_{\tilde{\tilde{A}}}(x_i)+\sum_{i=L+1}^N x_i \underline{\mu}_{\tilde{\tilde{A}}}(x_i)}{\sum_{i=1}^L \bar{\mu}_{\tilde{\tilde{A}}}(x_i)+\sum_{i=L+1}^N \underline{\mu}_{\tilde{\tilde{A}}}(x_i)},$$
$$y_r=\frac{\sum_{i=1}^R x_i \underline{\mu}_{\tilde{\tilde{A}}}(x_i)+\sum_{i=R+1}^N x_i \bar{\mu}_{\tilde{\tilde{A}}}(x_i)}{\sum_{i=1}^R \underline{\mu}_{\tilde{\tilde{A}}}(x_i)+\sum_{i=R+1}^N \bar{\mu}_{\tilde{\tilde{A}}}(x_i)}.$$
Here $L$ and $R$ are switch points which are calculated by KM
algorithm (\citealp{karn}; \citealp{mendL}). It is obvious that for
large $N$, i.e. for $(x_{i+1}-x_i)\rightarrow 0$, the discretization
of continuous domain is legitimate for computation of centroid.
Also, it is observed that (\citealp{mendL}) for IT2 FS with
symmetrical membership function, choice of $N$ has less effect on
computed centroid.
\subsection{Uncertainty bound (UB)}
\cite{wum1} provided inner- and outer-bound sets for type-reduced
set, which can not only be used to compute left and right end points
of the type-reduced set, but can also be used to derive the
defuzzified output of an IT2 FS. As compared to KM algorithm, this
method can performed without type-reduction and $x_i$s need not be
sorted, so that it removes computational burden of type-reduction.
The approximation of the type-reduced set by its inner- and
outer-bound sets is given by $[y_l,y_r]$, where
$y_l=(\underline{y}_l+\bar{y}_l)/2$ and
$y_r=(\underline{y}_r+\bar{y}_r)/2$, and corresponding defuzzified
output is $(y_l+y_r)/2$,
$$\bar{y}_l=\min\{y^{(0)},y^{(N)}\}$$
$$\underline{y}_r=\max\{y^{(0)},y^{(N)}\}$$
$$\underline{y}_l=\bar{y}_l-\left[\frac{\sum_{i=1}^N \left(\bar{\mu}_{\tilde{\tilde{A}}}(x_i)-\underline{\mu}_{\tilde{\tilde{A}}}(x_i)\right)}{\sum_{i=1}^N \bar{\mu}_{\tilde{\tilde{A}}}(x_i) \sum_{i=1}^N \underline{\mu}_{\tilde{\tilde{A}}}(x_i)}\times \frac{\sum_{i=1}^N \underline{\mu}_{\tilde{\tilde{A}}}(x_i) (x_i-x_1) \sum_{i=1}^N \bar{\mu}_{\tilde{\tilde{A}}}(x_i) (x_N-x_i)}{\sum_{i=1}^N \underline{\mu}_{\tilde{\tilde{A}}}(x_i) (x_i-x_1)+ \sum_{i=1}^N \bar{\mu}_{\tilde{\tilde{A}}}(x_i) (x_N-x_i)}\right]$$
$$\bar{y}_r=\underline{y}_r+\left[\frac{\sum_{i=1}^N \left(\bar{\mu}_{\tilde{\tilde{A}}}(x_i)-\underline{\mu}_{\tilde{\tilde{A}}}(x_i)\right)}{\sum_{i=1}^N \bar{\mu}_{\tilde{\tilde{A}}}(x_i) \sum_{i=1}^N \underline{\mu}_{\tilde{\tilde{A}}}(x_i)}\times \frac{\sum_{i=1}^N \bar{\mu}_{\tilde{\tilde{A}}}(x_i) (x_i-x_1) \sum_{i=1}^N \underline{\mu}_{\tilde{\tilde{A}}}(x_i) (x_N-x_i)}{\sum_{i=1}^N \bar{\mu}_{\tilde{\tilde{A}}}(x_i) (x_i-x_1)+ \sum_{i=1}^N \underline{\mu}_{\tilde{\tilde{A}}}(x_i) (x_N-x_i)}\right]$$
$$y^{(0)}=\frac{\sum_{i=1}^N x_i \underline{\mu}_{\tilde{\tilde{A}}}(x_i)}{\sum_{i=1}^N \underline{\mu}_{\tilde{\tilde{A}}}(x_i)}$$
$$y^{(N)}=\frac{\sum_{i=1}^N x_i \bar{\mu}_{\tilde{\tilde{A}}}(x_i)}{\sum_{i=1}^N \bar{\mu}_{\tilde{\tilde{A}}}(x_i)}.$$
\subsection{Nie-Tan (N-T) method}
\cite{nie} proposed a type-reduction method which is formulated
using the vertical-slice representation of an IT2 FS. In this
method, type reduction and defuzzification are performed together.
As of the previous two methods, if the domain of IT2 FS is
continuous then it is discreteized into finite number of points
$x_i$, $i=1,2,...,N$. Then the centroid (or defuzzified value) of
the IT2 FS $\tilde{\tilde{A}}$ can be expressed as
$$C(\tilde{\tilde{A}})=\frac{\sum_{i=1}^N x_i\bar{\mu}_{\tilde{\tilde{A}}}(x_i)+\sum_{i=1}^N x_i\underline{\mu}_{\tilde{\tilde{A}}}(x_i)}{\sum_{i=1}^N\bar{\mu}_{\tilde{\tilde{A}}}(x_i)+\sum_{i=1}^N\underline{\mu}_{\tilde{\tilde{A}}}(x_i)}.$$
The above formulation of the crisp output of an IT2 FS depends only
on the lower and upper bounds of its FOU. The computational
complexity of the N-T method is lower than the uncertainty bounds
method and the KM algorithm.
\subsection{Geometric centroid}
\cite{coup} introduced the idea of geometric centroid of an IT2 FS
by converting the region bounded by upper and lower membership
functions (which are piecewise linear) to a closed polygon. The
polygon consists of ordered coordinate points of the upper bound of
$\tilde{\tilde{A}}$ followed by the lower bound of
$\tilde{\tilde{A}}$ in reverse order. Let the polygon is given by
$((x_1,y_1),(x_2,y_2),...,(x_M,y_M))$, where $y_i$ is either
$\bar{\mu}_{\tilde{\tilde{A}}}(x_i)$ or
$\underline{\mu}_{\tilde{\tilde{A}}}(x_i)$ according to the position
of the coordinate point. Then the defuzzified output is taken as the
centroid (center of the polygon) of the polygon which is given by
$$\frac{\sum_{i=1}^{M-1} (x_i+x_{i+1})(x_i y_{i+1}-x_{i+1}y_i)}{3 \sum_{i=1}^{M-1}(x_i y_{i+1}-x_{i+1}y_i)}.$$

\section{Multi-objective optimization techniques}\label{secmulti}
The problem (\ref{morpt21})-(\ref{morpt23}) is a bi-objective
problem with one objective as maximization and another as
minimization. To solve this problem with defuzzified parameters we
apply different multi-objective optimization techniques, namely,
global criterion method, weighted sum method, desirability function
approach, fuzzy programming technique and NIMBUS which are discussed
briefly in this section. Consider a general multi-objective
optimization problem with some objectives to be maximized and some
others to be minimized:
\begin{eqnarray}&&\text{Max}~\{f_1(x),f_2(x),...,f_K(x)\}\label{eqmo1}\\&& \text{Min}~\{g_1(x),g_2(x),...,g_M(x)\}\label{eqmo2}\\ &&\text{s.t.}~x \in
D,\label{eqmo3}
\end{eqnarray}
where $D$ is the set of feasible solutions.\\ We use the following
notations in describing the methods:
$f_{i}^{\max}=\text{Max}~f_i(x)$, $f_{i}^{\min}=\text{Min}~f_i(x)$,
$g_{j}^{\max}=\text{Max}~g_j(x)$, and
$g_{j}^{\min}=\text{Min}~g_j(x)$, $i=1,2,...,K$ and $j=1,2,...,M$
subject to $x \in D$ in each case. We also denote the optimal
solution of single objective problem (considering only one objective
$f_{i}$ or $g_j$ ignoring all other objectives) as $x_{f_{i}}^{*}$
and $x_{g_{j}}^{*}$ respectively for $i=1,2,...,K$ and
$j=1,2,...,M$. The ideal objective vector for the above problem is
$(f_{1}^{\max},f_{2}^{\max},...,f_{K}^{\max},
g_{1}^{\min},g_{2}^{\min},...,g_{M}^{\min})$.
\subsection{Global criteria method}
By the method of global criteria (\citealp{zele}; \citealp{miett}) a
compromise solution is achieved by minimizing the sum of the
differences between ideal objective value and the respective
objective function values in the set of feasible solution. The ideal
objective value may be taken as the minimum value of the objective
function for minimization problem, and maximum value for
maximization problem obtained as by solving the multi-objective
problem as a single objective problem, considering each objective
individually. The method may be described by the following steps for
solving the multi-objective problem (\ref{eqmo1})-(\ref{eqmo3}):
\\Step 1: Construct single objective problems by taking each objective function individually.\\
Step 2: For each single objective problem, determine the ideal
objective vector $(f_{1}^{\max},f_{2}^{\max},...,f_{K}^{\max},
g_{1}^{\min},g_{2}^{\min},...,g_{M}^{\min})$
and corresponding values of \allowdisplaybreaks{$(f_{1}^{\min},f_{2}^{\min},...,f_{K}^{\min},g_{1}^{\max},\\g_{2}^{\max},...,g_{M}^{\max})$}.\\
Step 3: Formulate the following auxiliary problem using normalized
Minkowski distance ($L_p$ norm):
$$Min~G(x)$$ \hspace{6.5cm} s.t. $x \in D,$
$$G(x)=\text{Min}\left\{\sum_{i=1}^K \left(\frac{f_{i}^{\max}-f_{i}(x)}{f_{i}^{\max}-f_{i}^{\min}}\right)^p + \sum_{j=1}^M \left(\frac{g_{j}(x)-g_{j}^{\min}}{g_{j}^{\max}-g_{j}^{\min}}\right)^p\right\}^{\frac{1}{p}},$$
$$or,~G(x)=\text{Min}\left\{\sum_{i=1}^K \left(\frac{f_{i}^{\max}-f_{i}(x)}{f_{i}^{\max}}\right)^p + \sum_{j=1}^M \left(\frac{g_{j}(x)-g_{j}^{\min}}{g_{j}^{\min}}\right)^p\right\}^{\frac{1}{p}},$$
where $1\leq p\leq \infty$. An usual value of $p$ is 2.
\subsection{Weighted sum method}
In weighted sum method, multiple objectives are aggregated to
convert into a single objective by employing weight to each
objective. The weighting coefficients denote the relative importance
of the objectives. Now the values of different objectives may have
different order of magnitude, so it is necessary to normalize the
objectives, in order to make solution consistent with the weights as
assigned to the objectives. The objective functions may be converted
to their normal forms as follows:
$$f_{i}^{norm}=\frac{f_{i}(x)-f_{i}^{\min}}{f_{i}^{\max}-f_{i}^{\min}},~\text{(for objectives to be maximized)}$$
$$g_{j}^{norm}=\frac{g_{j}^{\max}-g_{j}(x)}{g_{j}^{\max}-g_{j}^{\min}},~\text{(for objectives to be minimized)}.$$
A weight $w_i$ is taken for every objective and then aggregated to
form the following problem:
$$\text{Max}~\left(\sum_{i=1}^K w_i f_{i}^{norm}+ \sum_{j=1}^M w_{K+j}~ g_{j}^{norm}\right)$$
\hspace{4.5cm} s.t. $x \in D,~w_i>0,~i=1,2,...,K+M,~\sum_{i=1}^{K+M}
w_i=1.$
\subsection{Desirability function approach}
By the desirability function approach (\citealp{akal};
\citealp{male}; \citealp{yeti}) each objective function $f_i$ is
transformed to a scale free desirability value $d_i$ $(0\leq d_i\leq
1)$ where $d_i=0$ represents completely undesirable response and
$d_i=1$ represents completely desirable or ideal response. Then
individual desirability values are aggregated into a single global
desirability index through a weighted geometric mean.\\ For the
objective function to be maximized its individual desirability
function $(d_i)$ is defined by
$$d_i=\left\{
        \begin{array}{ll}
          0, & \hbox{if $f_i< (f_i)_{wt}$;} \\
          \left(\frac{f_{i}(x)-(f_{i})_{wt}}{(f_{i})_{bt}-(f_{i})_{wt}}\right)^{k_i}, & \hbox{if $(f_i)_{wt}\leq f_i \leq (f_i)_{bt}$;} \\
          1, & \hbox{if $f_i > (f_i)_{bt}$}
        \end{array}
      \right.$$
where $(f_i)_{wt}$ and $(f_{i})_{bt}$ are the minimum (worst) and
the maximum (best) acceptable values of $f_i$, respectively. Here,
$k_i>0$ is the user-specified exponential parameter that determines
the shape (convex for $k_i < 1$ and concave for $k_i> 1$) of
desirability function. When $k_i=1$, the desirability function
increases linearly. Now for the objective function to be minimized
the individual desirability function $(s_j)$ is defined by
$$s_j=\left\{
        \begin{array}{ll}
          1, & \hbox{if $g_j< (g_j)_{bt}$;} \\
          \left(\frac{(g_j)_{wt}-g_j(x)}{(g_j)_{wt}-(g_j)_{bt}}\right)^{l_j}, & \hbox{if $(g_j)_{bt}\leq g_j \leq (g_j)_{wt}$;} \\
          0, & \hbox{if $g_j > (g_j)_{wt}$}
        \end{array}
      \right.$$
where $(g_j)_{wt}$ and $(g_j)_{bt}$ are the worst and the best
acceptable values of $g_j$, respectively, and $l_j>0$.\\ The overall
desirability $d$ which combines the individual desirability values
into a single response is defined as the weighted geometric mean of
all the individual desirability values:
$$d=\left(d_1^{w_1}\times \ldots \times d_K^{w_K}\times s_1^{w_{K+1}}\times \ldots \times s_M^{w_{K+M}}\right)^{1/(w_1+w_2+...+w_{K+M})},$$
where $w_r$ $(r=1,2,...,K+M)$ represents relative importance
(\citealp{akal}) that varies from the least important a value of 1,
to the most important a value of 5. The overall desirability
$d~(0\leq d\leq 1)$ has to be maximized subject to the constraints
of the problem to find the most desirable solution.\\ \textbf{Note:}
It is obvious that maximum (best) acceptable value for an objective
should be its optimal value as obtained by solving the problem as
single objective, e.g. $(f_{i})_{bt}=f_i^{\max}$ and
$(g_{j})_{bt}=g_j^{\min}$. We propose to take minimum (worst)
acceptable value for an objective to be maximized as the minimum of
the values of that objective function evaluated at the optimal
solutions of all the single objective problems, i.e.
$$(f_{i})_{wt}=\text{Min}\{f_{i}(x_{f_{1}}^{*}),...,f_{i}(x_{f_{K}}^{*}),f_{i}(x_{g_{1}}^{*}),...,f_{i}(x_{g_{M}}^{*})\}$$
and for an objective to be minimized
$$(g_{j})_{wt}=\text{Max}\{g_{j}(x_{f_{1}}^{*}),...,g_{j}(x_{f_{K}}^{*}),g_{j}(x_{g_{1}}^{*}),...,g_{j}(x_{g_{M}}^{*})\}.$$
\subsection{Fuzzy programming technique} \cite{zimm} (see also \citealp{bit}; \citealp{kund1}) introduced fuzzy linear programming approach to solve
multi-objective problem, and he showed that fuzzy linear programming
always gives efficient solutions and an optimal compromise solution.
This method consists of the following steps to solve the
multi-objective problem (\ref{eqmo1})-(\ref{eqmo3}):
\\Step 1: Solve the problem taking each objective individually (ignoring all other
objectives) and obtain the corresponding optimal solutions
$x_{f_{i}}^{*}$, $i=1,2,...,K$ and $x_{g_{j}}^{*}$, $j=1,2,...,M$.\\
Step 2: Calculate the values of each objective function at all these
optimal solutions $x_{f_{i}}^{*}$ and $x_{g_{j}}^{*}$ and find the
upper and lower bound for each objective given by\\
$U_{f_{i}}=f_{i}(x_{f_{i}}^{*})$ and $L_{f_{i}}=\text{Min}\{f_{i}(x_{f_{1}}^{*}),...,f_{i}(x_{f_{K}}^{*}),f_{i}(x_{g_{1}}^{*}),...,f_{i}(x_{g_{M}}^{*})\}$,\\
$U_{g_{j}}=\text{Max}\{g_{j}(x_{f_{1}}^{*}),...,g_{j}(x_{f_{K}}^{*}),g_{j}(x_{g_{1}}^{*}),...,g_{j}(x_{g_{M}}^{*})\}$ and $L_{g_{j}}=g_{j}(x_{g_{j}}^{*})$, respectively.\\
Step 3: Construct the linear membership functions corresponding to
each objective as
$$\mu_{f_{i}}(f_{i})=\left\{
                                      \begin{array}{ll}
                                        0, & \hbox{if $f_{i}\leq L_{f_{i}}$;} \\
                                        \frac{f_{i}(x)-L_{f_{i}}}{U_{f_{i}}-L_{f_{i}}}, & \hbox{if $L_{f_{i}}\leq f_{i}\leq U_{f_{i}}$;} \\
                                        1, & \hbox{if $f_{i}\geq U_{f_{i}}$,}
                                      \end{array}
                                    \right.$$
$$\mu_{g_{j}}(g_{j})=\left\{
                                      \begin{array}{ll}
                                        1, & \hbox{if $g_{j}\leq L_{g_{j}}$;} \\
                                        \frac{U_{g_{j}}-g_{j}(x)}{U_{g_{j}}-L_{g_{j}}}, & \hbox{if $L_{g_{j}}\leq g_{j}\leq U_{g_{j}}$;} \\
                                        0, & \hbox{if $g_{j}\geq U_{g_{j}}$.}
                                      \end{array}
                                    \right.$$
Step 4: Formulate fuzzy linear programming problem using max-min
operator for the multi-objective problem as
\begin{eqnarray*}&&\text{Max}~\min_{i,j}\{\mu_{f_{i}}(f_{i}),\mu_{g_{j}}(g_{j})\}\\
&&s.t.~x \in D,
\end{eqnarray*}
i.e.
$$Max ~\lambda$$ $$\text{subject~to}~
\mu_{f_{i}}(f_{i})=(f_{i}(x)-L_{f_{i}})/(U_{f_{i}}-L_{f_{i}})\geq
\lambda,~i=1,...,K,$$
$$\mu_{g_{j}}(g_{j})=(U_{g_{j}}-g_{j}(x))/(U_{g_{j}}-L_{g_{j}})\geq \lambda,~j=1,...,M,$$ $$x \in
D,~ \lambda\geq 0.$$\\ Step 5: Solve the reduced problem of step 4
by a linear optimization technique, and the optimum compromise
solutions are obtained.
\subsection{NIMBUS}
\cite{miet} introduced a methodology known as NIMBUS method for
solving interactive multi-objective optimization problems. The
solution process is based on the classification of objective
functions. In this method, several scalarizing functions are
formulated based on the objective functions and the preference
information specified by the decision maker, and they usually
generate Pareto optimal (PO) solutions for the original problem. In
classification, first objective function values are calculated at
the current PO decision vector, say $x^c$, and then every objective
function is put into one of the classes based on desirable changes
in the objective function values. There are five different classes
for each of the objective functions $g_i$ (say) whose values -
should be improved as much as possible $(i\in I^{imp})$, should be
improved till some desired aspiration level $\bar{g}_i< g_{i}(x^c)$
(for minimization problem) $(i\in I^{asp})$, is satisfactory at the
moment $(i\in I^{sat})$, is allowed to get worse until a value
$\varepsilon_i$ $(i\in I^{bound})$, and can change freely at the
moment $(i\in I^{free})$. A classification is feasible only if
$$I^{imp} \cup I^{asp}\neq \emptyset ~ \text{and}~ I^{bound} \cup I^{free}\neq \emptyset.$$
A scalarized subproblem is then formed based on the classification
and the corresponding aspiration levels and upper bounds as follows
(for minimization problem):
\begin{eqnarray*}&&\text{Min}~\text{Max}_{i\in I^{imp}, j\in
I^{asp}}~\left[\frac{g_i(x)-g_i(x_{g_i}^*)}{g_{i}^{nad}-g_i(x_{g_i}^*)},
\frac{g_j(x)-\bar{g}_j}{g_{j}^{nad}-g_j(x_{g_j}^*)}\right]+\rho
\sum_{i=1}^M \frac{g_i(x)}{g_{i}^{nad}-g_i(x_{g_i}^*)}\\ &&
\text{s.t.}~g_i(x)\leq g_{i}(x^c)~\forall~i\in I^{imp}\cup
I^{asp}\cup I^{sat},\\&& g_i(x)\leq~\varepsilon_i~\forall i\in
I^{bound},
\\&& x \in D,
\end{eqnarray*}
where $\rho>0$ is an augmentation coefficient and is relatively a
small scalar. Solution of the scalarized problem is either weakly PO
or PO according to the augmentation coefficient is used or not used.
\cite{miet} implemented NIMBUS method as a WWW-NIMBUS software system which is accessible at http://nimbus.mit.jyu.fi/.\\
Convergence indicator: To discuss the convergence of the
multi-objective optimization procedure or to measure the quality of
the solution, we adopt a convergence indicator or measure of
performances, namely Convergence Metric or Distance Metric $d$ to
find Euclidean distance (normalized) between ideal solution and
compromise solution. This indicator will measure closeness of the
obtained compromise objective values with the respective ideal
objective values. The smaller this metric value, the better is the
convergence towards the ideal solution.
\section{Numerical Experiment}
To illustrate the MORRAP (\ref{morpt21})-(\ref{morpt23}), i.e. the
problem (\ref{eqmorp1})-(\ref{eqmorp5}) with imprecise component
reliabilities represented as IT2 FNs, we consider a
reliability-redundancy allocation problem on a pharmaceutical plant
(for details see Garg and Sharma, 2013), where two objectives are
maximization of system reliability and minimization of system cost.
The mathematical formulation of the bi-objective problem is given by
(\ref{morpex1})-(\ref{morpex5}) with the input parameters given in
Table 1.
\begin{eqnarray}\label{morpex1}&&\text{Max}~\tilde{\tilde{R}}_{s}=\prod_{i=1}^{10}~\left[1-(1-\tilde{\tilde{r}}_i)^{n_i}\right],\\
\label{morpex2}&&Min~\tilde{\tilde{C}}_s=\sum_{i=1}^{10}
\alpha_{i}\left(\frac{-T}{ln(\tilde{\tilde{r}}_i)}\right)^{\beta_{i}}(n_i+\exp(n_i/4)),\end{eqnarray}
\begin{eqnarray}\label{morpex3} \text{subject to}
&&V_{s}=\sum_{i=1}^{10} v_{i} n_{i}^2~\leq V,\\
\label{morpex4}&&W_{s}=\sum_{i=1}^{10} w_{i}(n_i\cdot
\exp(n_i/4))\leq W,
\end{eqnarray}
\begin{equation}\label{morpex5}1\leq n_i \leq 3, n_i\in \mathbb{Z}^+, i=1,2,...,10,\end{equation}
where $\tilde{\tilde{r}}_i$ is represented by IT2 FN having support
$\subseteq [0.5, 1-10^{-6}]$. The IT2 FN $\tilde{\tilde{r}}_i$,
$i=1,2,...,10$ are generated using the Algorithm presented in
Section \ref{secMOIT2} and are given in Table 2, where approximate
reliabilities are given by `about $r_i$', $i=1,2,...,10$, and
$r_1=0.55$, $r_2=0.60$, $r_3=0.65$, $r_4=0.70$, $r_5=0.75$,
$r_6=0.80$, $r_7=0.85$, $r_8=0.90$, $r_9=0.92$, $r_{10}=0.95$. We
apply various type-reduction strategies and defuzzification
techniques as discussed in Section \ref{sectype} to obtain
corresponding defuzzified values of IT2 FNs and are presented in
Table 3. In that table, for the KM algorithm and the UB method we
also provide left and right end points of the centroid and
uncertainty bounds respectively, along with the corresponding
defuzzified values. From the defuzzified values in Table 3, it is
observed that KM algorithm, uncertainty bound approach and N-T
method give more similar result as compared to the geometric
centroid approach.\par
\begin{table}[h]
\begin{center} {Table 1: Input parameters \\
\begin{tabular}{cccccccc}
\hline Components & $10^5 \alpha_{i}$  & $\beta_{i}$ & $v_{i}$ &
$w_{i}$ & $V$ & $W$ & $T$\\\hline
1& 0.611360 & 1.5 & 4.0 & 9.0 &  &  & \\
2& 4.032464 & 1.5 & 5.0 & 7.0 &  &  &\\
3& 3.578225 & 1.5 & 3.0 & 5.0 &  &  &\\
4& 3.654303 & 1.5 & 2.0 & 9.0 &  &  &\\
5& 1.163718 & 1.5 & 3.0 & 9.0 & 289 & 483 & 1000\\
6& 2.966955 & 1.5 & 4.0 & 10.0 &  &  &\\
7& 2.045865 & 1.5 & 1.0 & 6.0 &  &  &\\
8& 2.649522 & 1.5 & 1.0 & 5.0 &  &  &\\
9& 1.982908 & 1.5 & 4.0 & 8.0 &  &  &\\
10& 3.516724 & 1.5 & 4.0 & 6.0 &  &  &\\
\hline
\end{tabular} }
\end{center}
\end{table}
\begin{table}[h]
\begin{center} {Table 2: IT2 FN $\tilde{\tilde{r}}_i$\\
\begin{tabular}{|c|c|}
\hline $\tilde{\tilde{r}}_1$ & ((0.511813,0.55,0.893671),(0.542672,0.55,0.615958))\\
$\tilde{\tilde{r}}_2$ & ((0.523627,0.60,0.905484),(0.585344,0.60,0.658620))\\
$\tilde{\tilde{r}}_3$ & ((0.535440,0.65,0.917298),(0.628017,0.65,0.701292))\\
$\tilde{\tilde{r}}_4$ & ((0.547254,0.70,0.929111),(0.670689,0.70,0.743965))\\
$\tilde{\tilde{r}}_5$ & ((0.559067,0.75,0.940925),(0.713361,0.75,0.786637))\\
$\tilde{\tilde{r}}_6$ & ((0.570880,0.80,0.952738),(0.756034,0.80,0.829309))\\
$\tilde{\tilde{r}}_7$ & ((0.582694,0.85,0.964552),(0.798706,0.85,0.871981))\\
$\tilde{\tilde{r}}_8$ & ((0.594508,0.90,0.976365),(0.841378,0.90,0.914654))\\
$\tilde{\tilde{r}}_9$ & ((0.599233,0.92,0.981091),(0.858447,0.92,0.931723))\\
$\tilde{\tilde{r}}_{10}$ & ((0.606321,0.95,0.988170),(0.884050,0.95,0.957326))\\
\hline
\end{tabular} }
\end{center}
\end{table}
\begin{table}[h]
\begin{center} {Table 3: Defuzzified values with different type-reduction strategies \\
\begin{tabular}{|c|cccc|}
\hline IT2 FN  & Centroid value & Uncertainty bound & Defuzzified value & Geometric centroid \\
               & (KM Algorithm) &  & (N-T method) &  \\\hline
$\tilde{\tilde{r}}_1$ & [0.559313,0.685104] & [0.54701,0.741079] &  &  \\
 & 0.622208 & 0.644044 & 0.638117 & 0.671368 \\
$\tilde{\tilde{r}}_2$  & [0.594175,0.714798] & [0.584012,0.761516] &  &  \\
  & 0.654486 & 0.672764 & 0.666158 & 0.691025 \\
$\tilde{\tilde{r}}_3$  & [0.628406,0.744975] & [0.614688,0.780418] &  &  \\
 & 0.686690 & 0.697553 & 0.694166 & 0.710682 \\
$\tilde{\tilde{r}}_4$  & [0.661416,0.775753] & [0.649731,0.798093] &  &  \\
 & 0.718584 & 0.723912 & 0.722142 & 0.730339 \\
$\tilde{\tilde{r}}_5$  & [0.693230,0.806764] & [0.685508,0.814486] &  &  \\
 & 0.749997 & 0.749997 & 0.749997 & 0.749996 \\
$\tilde{\tilde{r}}_6$  & [0.724241,0.838579] & [0.701899,0.850265] &  &  \\
 & 0.781410 & 0.776082 & 0.777853 & 0.769654 \\
$\tilde{\tilde{r}}_7$  & [0.755019,0.87159] & [0.719574,0.885308] &  &  \\
 & 0.813304 & 0.802441 & 0.805828 & 0.789311 \\
$\tilde{\tilde{r}}_8$  & [0.785194,0.905821] & [0.738475,0.919584] &  &  \\
 & 0.845507 & 0.829029 & 0.833836 & 0.808968 \\
$\tilde{\tilde{r}}_9$  & [0.795185,0.919755] & [0.744763,0.932876] &  &  \\
 & 0.857470 & 0.838819 & 0.844481 & 0.816831 \\
$\tilde{\tilde{r}}_{10}$  & [0.814883,0.940682] & [0.758908,0.952984] &  &  \\
 & 0.877782 & 0.855946 & 0.861875 & 0.828625 \\
\hline
\end{tabular} }
\end{center}
\end{table}
With the defuzzified values as given in Table 3, we solve the
bi-objective problem (\ref{morpex1})-(\ref{morpex5}) by applying
different multi-objective techniques as discussed in Section
\ref{secmulti}. The results are obtained using standard optimization
solver LINGO which is based on gradient based optimization -
Generalized Reduced Gradient (GRG) technique. Tables 4-7 provide the
solution of the problem with five different multi-objective
techniques where the defuzzified values are obtained by KM
algorithm, uncertainty bound, N-T method and geometric centroid,
respectively. From the results (Tables 4-7) it is observed that, the
subsystem comprising of components with lower reliability (e.g.
subsystem 1) is associated higher redundancy to increase the
reliability of the system. To the contrary, the subsystem comprising
of components with higher reliability (e.g. subsystem 10) is
associated fewer redundancy to reduce the cost of the system. Also,
the two objectives of the problem are conflicting to each other, so
we can only derived compromise solutions (as seen from the results
presented in Tables 4-7). For multi-objective problem with
conflicting objectives it is not easy to compare the results as
obtained by different methods. However, different results in hand
gives more flexibility to a decision maker (DM) to choose
appropriate result according to his/her preference or as situation
demand. For instance if DM emphasizes more preference on reliability
maximization over cost minimization, then DM may consider the
results obtained by desirability function approach and weighted sum
method. If DM's preference is more on cost minimization, then the
results obtained by fuzzy programming approach and NIMBUS can be
chosen. Whereas, if DM's determination not to give preference to one
objective over the other, then the results obtained by global
criteria method in $L_2$ norm can be chosen. One can also measure
the quality of the solution, by adopting a convergence indicator or
measure of performances. Here we choose Convergence Metric or
Distance Metric to find Euclidean distance (normalized) between
compromise objective values and the respective ideal objective
values. The smaller this metric value, the better is the convergence
towards the ideal solution. For the solutions obtained by the
different multi-objective optimization techniques as presented in
Table 4, the values of the corresponding normalized Euclidean
distances are calculated as 0.6075097, 0.7629309, 0.9145547,
0.7629310, 0.5541247, and 0.5609847 respectively. Similar
observations can be made for the solutions obtained by the different
multi-objective optimization techniques as presented in Tables 5, 6
and 7.
\begin{table}[h]
\begin{center} {Table 4: Solution with different multi-objective
optimization techniques for the problem with defuzzified values
obtained using KM
Algorithm\\
\begin{tabular}{|c|l|}
\hline Individual optimal value & Max $R_s=0.8317749$, Min $C_s=181.2395$\\
\hline Multi-objective  Method & \hspace{3.5cm} Compromise solution \\
\hline Global criteria ($p=2$) & $R_s=0.6846485$, $C_s=286.5739$,
$n_1=5$,
$n_2=3$, $n_3=3$, $n_4=3$, \\
& $n_5=3$, $n_6=2$, $n_7=2$, $n_8=2$, $n_9=2$, $n_{10}=1$.\\\hline
Weighted sum & $R_s=0.7683246$, $C_s=318.8198$, $n_1=5$,
$n_2=3$, $n_3=3$, $n_4=3$, \\
(with equal weights) &$n_5=3$, $n_6=2$, $n_7=2$, $n_8=2$, $n_9=2$,
$n_{10}=2$.\\\hline Desirability function & $R_s=0.829084$,
$C_s=346.9919$, $n_1=4$, $n_2=3$, $n_3=4$, $n_4=3$,\\
$(t_1=1,t_2=0.1,w_1=w_2)$ & $n_5=2$, $n_6=2$, $n_7=2$, $n_8=1$, $n_9=2$, $n_{10}=1$.\\
$(t_1=0.5,t_2=0.1,w_1=w_2)$ & $R_s=0.768324$, $C_s=318.8198$,
$n_1=5$, $n_2=3$, $n_3=3$, $n_4=3$,\\
& $n_5=3$, $n_6=3$, $n_7=3$, $n_8=2$, $n_9=2$, $n_{10}=2$.\\ \hline
Fuzzy programming & $R_s=0.5319160$, $C_s=257.5089$, $n_1=5$,
$n_2=3$, $n_3=3$, $n_4=2$,\\
& $n_5=3$, $n_6=2$, $n_7=2$, $n_8=2$, $n_9=2$, $n_{10}=2$.\\ \hline
NIMBUS & $R_s=0.5306198$, $C_s=258.901$,
$n_1=4$, $n_2=3$, $n_3=3$, $n_4=2$, \\
& $n_5=2$, $n_6=2$, $n_7=2$, $n_8=2$, $n_9=1$, $n_{10}=1$.\\\hline
\end{tabular} }
\end{center}
\end{table}
\begin{table}[h]
\begin{center} {Table 5: Solution with different multi-objective
optimization techniques for the problem with defuzzified values
obtained using Uncertainty bound\\
\begin{tabular}{|c|l|}
\hline Individual optimal value & Max $R_s=0.8382419$, Min $C_s=160.4723$\\
\hline Multi-objective  Method & \hspace{3.5cm} Compromise solution \\
\hline Global criteria ($p=2$) & $R_s=0.6641386$, $C_s=262.7524$,
$n_1=5$,
$n_2=3$, $n_3=3$, $n_4=3$,\\
& $n_5=3$, $n_6=2$, $n_7=2$, $n_8=2$, $n_9=2$, $n_{10}=1$.\\\hline
Weighted sum & $R_s=0.7598104$, $C_s=287.4911$,
$n_1=5$, $n_2=3$, $n_3=3$, $n_4=3$,\\
(with equal weights) & $n_5=3$, $n_6=2$, $n_7=2$, $n_8=2$, $n_9=2$,
$n_{10}=2$.\\\hline Desirability function & $R_s=0.8082213$,
$C_s=306.3102$, $n_1=4$, $n_2=3$, $n_3=3$, $n_4=3$,\\
$(t_1=1,t_2=0.1,w_1=w_2)$ & $n_5=3$, $n_6=3$, $n_7=3$, $n_8=2$, $n_9=2$, $n_{10}=2$.\\
$(t_1=0.5,t_2=0.1,w_1=w_2)$ & $R_s=0.7598104$, $C_s=287.4911$,
$n_1=5$, $n_2=3$, $n_3=3$, $n_4=3$,\\
& $n_5=3$, $n_6=2$, $n_7=2$, $n_8=2$, $n_9=2$, $n_{10}=2$.\\\hline
Fuzzy programming & $R_s=0.5160557$, $C_s=234.8222$, $n_1=5$,
$n_2=2$, $n_3=2$, $n_4=2$,
\\ & $n_5=2$, $n_6=2$, $n_7=2$, $n_8=2$, $n_9=2$, $n_{10}=1$.\\\hline
NIMBUS & $R_s=0.5160557$,
$C_s=234.8222$, $n_1=5$, $n_2=2$, $n_3=2$, $n_4=2$,\\
& $n_5=2$, $n_6=2$, $n_7=2$, $n_8=2$, $n_9=2$, $n_{10}=1$.\\ \hline
\end{tabular} }
\end{center}
\end{table}
\begin{table}[h]
\begin{center} {Table 6: Solution with different multi-objective
optimization techniques for the problem with defuzzified values
obtained using N-T Method\\
\begin{tabular}{|c|l|}
\hline Individual optimal value & Max $R_s=0.8363644$, Min $C_s=165.4758$\\
\hline Multi-objective Method & \hspace{3.5cm} Compromise solution \\
\hline Global criteria ($p=2$) & $R_s=0.6698056$, $C_s=268.3749$,
$n_1=5$,
$n_2=3$, $n_3=3$, $n_4=3$,\\
& $n_5=3$, $n_6=2$, $n_7=2$, $n_8=2$, $n_9=2$, $n_{10}=1$.\\\hline
Weighted sum & $R_s=0.7623225$, $C_s=294.8568$, $n_1=5$,
$n_2=3$, $n_3=3$, $n_4=3$,\\
(with equal weights) & $n_5=3$, $n_6=2$, $n_7=2$, $n_8=2$, $n_9=2$,
$n_{10}=2$.\\\hline Desirability function & $R_s=0.8091350$,
$C_s=314.1297$, $n_1=4$, $n_2=3$, $n_3=3$, $n_4=3$,\\
$(t_1=1,t_2=0.1,w_1=w_2)$& $n_5=3$, $n_6=3$, $n_7=3$, $n_8=2$, $n_9=2$, $n_{10}=2$.\\
$(t_1=0.5,t_2=0.1,w_1=w_2)$ & $R_s=0.7623225$, $C_s=294.8568$,
$n_1=5$, $n_2=3$, $n_3=3$, $n_4=3$,\\
&$n_5=3$, $n_6=2$, $n_7=2$, $n_8=2$, $n_9=2$, $n_{10}=2$.\\\hline
Fuzzy programming & $R_s=0.5180679$, $C_s=240.9737$, $n_1=5$,
$n_2=2$, $n_3=2$,
$n_4=2$,\\
& $n_5=2$, $n_6=2$, $n_7=2$, $n_8=2$, $n_9=2$, $n_{10}=1$.\\\hline
NIMBUS & $R_s=0.5191948$,
$C_s=242.3131$, $n_1=3$, $n_2=3$, $n_3=3$, $n_4=2$,\\
& $n_5=3$, $n_6=2$, $n_7=2$, $n_8=1$, $n_9=2$, $n_{10}=1$.\\ \hline
\end{tabular} }
\end{center}
\end{table}
\begin{table}[h]
\begin{center} {Table 7: Solution with different multi-objective
optimization techniques for the problem with defuzzified values
obtained using geometric centroid\\
\begin{tabular}{|c|l|}
\hline Individual optimal value & Max $R_s=0.8470077$, Min $C_s=143.4406$\\
\hline Multi-objective Method & \hspace{3.5cm} Compromise solution \\
\hline Global criteria ($p=2$) & $R_s=0.6561468$, $C_s=243.3404$,
$n_1=4$,
$n_2=3$, $n_3=2$, $n_4=2$,\\
& $n_5=3$, $n_6=2$, $n_7=2$, $n_8=2$, $n_9=2$, $n_{10}=2$.\\\hline
Weighted sum & $R_s=0.7446174$, $C_s=262.6584$, $n_1=5$,
$n_2=3$, $n_3=3$, $n_4=3$,\\
(with equal weights) & $n_5=3$, $n_6=2$, $n_7=2$, $n_8=2$, $n_9=2$,
$n_{10}=2$. \\\hline Desirability function & $R_s=0.8215322$,
$C_s=289.9504$, $n_1=4$, $n_2=3$, $n_3=3$, $n_4=3$,\\
$(t_1=1,t_2=0.1,w_1=w_2)$ & $n_5=3$, $n_6=3$, $n_7=3$, $n_8=2$, $n_9=3$, $n_{10}=2$.\\
$(t_1=0.5,t_2=0.1,w_1=w_2)$ & $R_s=0.7719188$, $C_s=270.9126$,
$n_1=5$, $n_2=3$, $n_3=3$, $n_4=3$,\\
&$n_5=3$, $n_6=2$, $n_7=3$, $n_8=2$, $n_9=2$, $n_{10}=2$.\\\hline
Fuzzy programming & $R_s=0.5220752$, $C_s=216.3870$, $n_1=4$,
$n_2=2$, $n_3=2$, $n_4=2$,\\ & $n_5=$, $n_6=3$, $n_7=2$, $n_8=2$,
$n_9=2$, $n_{10}=1$.
\\\hline
NIMBUS & $R_s=0.5008404$,
$C_s=221.3302$, $n_1=4$, $n_2=3$, $n_3=3$, $n_4=2$,\\
& $n_5=3$, $n_6=2$, $n_7=2$, $n_8=1$, $n_9=2$, $n_{10}=1$.\\ \hline
\end{tabular} }
\end{center}
\end{table}
\begin{figure}\label{figrelr}\begin{center}
  \includegraphics[width=10cm]{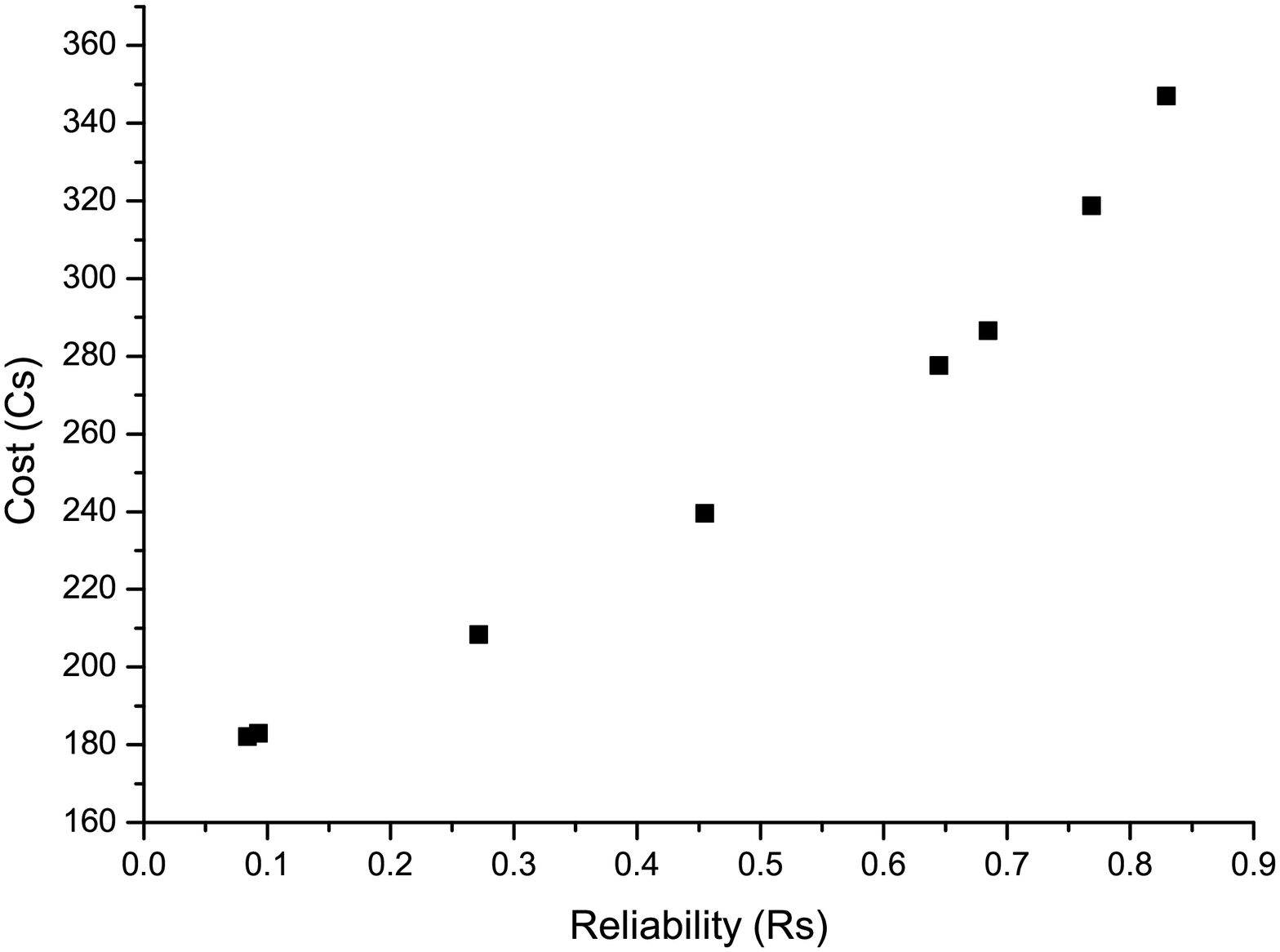}\\\vspace{-0.6cm}
  \caption{Nondominated solutions.}\end{center}
\end{figure}

In Tables 4-7, we have given single solution for different
multi-objective optimization methods by adopting suitable criteria,
e.g. for weighted sum method we chose equal weights for each of the
objectives; for global criterion method we use $L_2$ norm, etc.
However, Pareto optimality can be observed if one wishes to do so.
Here we construct a Pareto front (non-dominated solutions)
considering the weighted sum approach by assigning different
weights, i.e., $w_1$ and $w_2$ for the objectives $R_s$ and $C_s$
respectively, where $w_1,w_2\in[0,1]$ and $w_1+w_2=1$. The Pareto
front is depicted in Fig. 3.
\par Next we solve the
problem (\ref{eqmorp1})-(\ref{eqmorp5}) with the component
reliabilities represented as T1 FNs having support $\subseteq [0.5,
1-10^{-6}]$, instead of IT2 FNs. The T1 FNs $\tilde{r}_i$,
$i=1,2,...,10$ can be generated using the Steps 1-3 of the Algorithm
presented in Section \ref{secMOIT2}. Our intensity is to compare the
results of MORRAP with uncertain component reliabilities represented
as IT2 FNs and that of same problem if one represents uncertain
component reliabilities by T1 FNs. For this purpose, in Table 8, we
present the solution of MORRAP with T1 FNs where defuzzified values
are obtained using centroid value of T1 FN. It is to be noted that
the centroid of a T1 FN $(r^l,r,r^u)$ is given by $(r^l+r+r^u)/3$.
For comparison, in the Table 8, we also display the solution of the
problem with IT2 FNs where defuzzified (centroid) values are
obtained using KM Algorithm. To avoid biasedness in the comparative
study we obtain the solutions using five different multi-objective
optimization techniques. The results are also display in the Fig. 4
for better realization. From the Table 8 and Fig. 4, it is observed
that modeling uncertain parameters (reliabilities) using IT2 FNs
leads to the better performance than that of using T1 FNs, i.e. we
can model system with higher system reliability and less system
cost. It is to be noted here that for the result obtained using
global criteria method, system reliability for the problem with IT2
FNs is slightly lower than that of with T1 FNs, but in this case
system cost is also much lower than the problem with T1 FNs.
\begin{figure}\begin{center}
  \includegraphics[width=15.5cm]{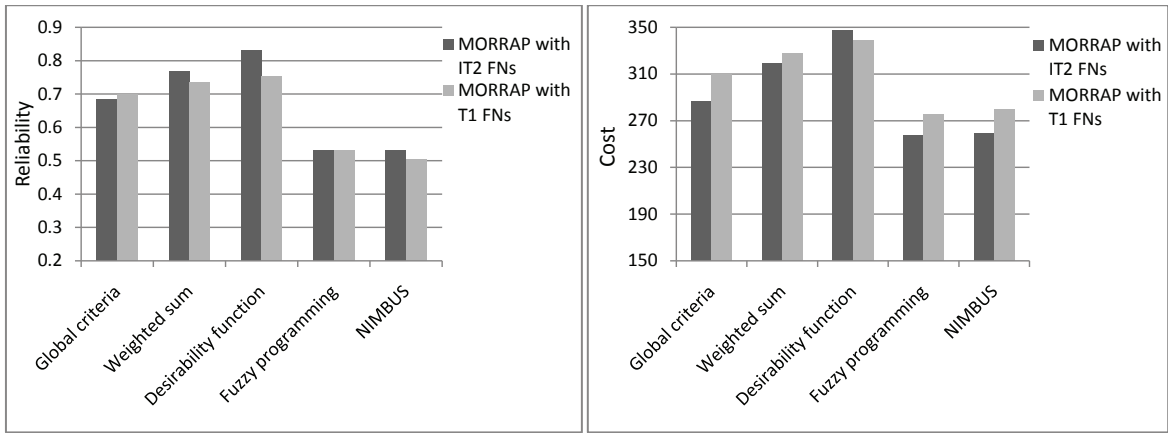}\\\caption{Comparative result of MORRAP with IT2 FNs and T1 FNs}
 \end{center}
\end{figure}

\begin{table}[h]
\begin{center} {Table 8: Solution of MORRAP with IT2 FNs and T1 FNs\\
\resizebox{\textwidth}{!}{%
\begin{tabular}{|c|l|l|}
\hline & MORRAP with IT2 FNs & MORRAP with T1 FNs\\
\hline Individual optimal value & Max $R_s=0.8317749$, & Max $R_s=0.8241383$, \\
& Min $C_s=181.2395$ & Min $C_s=203.9562$\\
\hline Multi-objective  Method & Compromise solution & Compromise solution \\
\hline Global criteria ($p=2$) & $R_s=0.6846485$, $C_s=286.5739$, & $R_s=0.6974577$, $C_s=309.8418$,\\
& $n_1=5$, $n_2=3$, $n_3=3$, $n_4=3$, $n_5=3$, & $n_1=5$, $n_2=3$, $n_3=3$, $n_4=3$, $n_5=3$,\\
& $n_6=2$, $n_7=2$, $n_8=2$, $n_9=2$, $n_{10}=1$. & $n_6=2$,
$n_7=2$, $n_8=2$, $n_9=2$, $n_{10}=1$.\\\hline Weighted sum &
$R_s = 0.7683246$, $C_s = 318.8198$, &  $R_s=0.7349505$, $C_s=327.4596$,\\
& $n_1=5$, $n_2=3$, $n_3=3$, $n_4=3$, $n_5=3$, & $n_1=4$, $n_2=4$,
$n_3=3$, $n_4=3$, $n_5=3$,\\
(with equal weights) & $n_6=2$, $n_7=2$, $n_8=2$, $n_9=2$, $n_{10}=2$. & $n_6=3$, $n_7=2$, $n_8=2$, $n_9=2$, $n_{10}=1$.\\
\hline Desirability function & $R_s=0.829084$,
$C_s=346.9919$, & $R_s=0.7542184$, $C_s=338.6593$,\\
$(t_1=1,t_2=0.1,$ & $n_1=4$, $n_2=3$, $n_3=4$, $n_4=3$, $n_5=3$, &
$n_1=4$, $n_2=4$, $n_3=3$, $n_4=3$, $n_5=3$,\\ $w_1=w_2)$ & $n_6=3$, $n_7=3$, $n_8=2$, $n_9=2$, $n_{10}=2$. & $n_6=3$, $n_7=3$, $n_8=2$, $n_9=2$, $n_{10}=1$.\\
\hline Fuzzy programming & $R_s=0.5319160$, $C_s=257.5089$, & $R_s=0.5321748$, $C_s=275.6192$,\\
& $n_1=5$, $n_2=3$, $n_3=3$, $n_4=2$, $n_5=2$, & $n_1=4$, $n_2=3$,
$n_3=2$, $n_4=2$, $n_5=3$,\\ & $n_6=2$, $n_7=2$, $n_8=1$, $n_9=2$,
$n_{10}=1$. & $n_6=2$, $n_7=2$, $n_8=2$, $n_9=1$, $n_{10}=1$.\\
\hline NIMBUS & $R_s=0.5306198$, $C_s=258.901$, & $R_s=0.5040566$,
$C_s=280.0497$,\\ & $n_1=4$, $n_2=3$, $n_3=3$, $n_4=2$,  $n_5=2$, & $n_1=4$, $n_2=3$, $n_3=3$, $n_4=3$, $n_5=4$,\\
& $n_6=2$, $n_7=2$, $n_8=2$, $n_9=1$, $n_{10}=1$. & $n_6=1$,
$n_7=2$, $n_8=2$, $n_9=1$, $n_{10}=1$.\\\hline
\end{tabular} }}
\end{center}
\end{table}

\section{Conclusion}
In this paper, we consider a multi-objective reliability-redundancy
allocation problem (MORRAP) of a series-parallel system. Here,
system reliability has to be maximized, and system cost has to be
minimized simultaneously subject to limits on weight, volume, and
redundancy level. Use of redundant components is commonly adapted
approach to increase reliability of a system. However, incorporation
of more redundant components may increase the cost of the system,
for which optimal redundancy is mainly concerned for the economical
design of system. Also, the component reliabilities in a system
cannot always be precisely measured as crisp values, but may be
determined as approximate values or approximate intervals with
imprecise endpoints. To deal with impreciseness, the presented
problem is formulated with the component reliabilities represented
as IT2 FNs which are more flexible and appropriate to model
impreciseness over usual or T1 FNs.

To solve MORRAP with interval type-2 fuzzy parameters, we first
apply various type-reduction and corresponding defuzzification
techniques, and obtain corresponding defuzzified values to observe
the effect of different type-reduction strategies. We illustrate the
problem with a real-world MORRAP on pharmaceutical plant. The
objectives of the problem are conflicting with each other, and so
one can obtain compromise solution in the sense that individual
optimal solution can not be reached together. To deal with this, we
apply five different multi-objective optimization techniques in the
view that different results in hand give more flexibility to a
decision maker to choose appropriate result according to his/her
preference or as situation demand. We also solve the MORRAP with the
uncertain (imprecise) component reliabilities represented as T1 FNs,
and observe that modeling impreciseness using IT2 FNs leads to
better performance than that of using T1 FNs. The present
investigation has been done by modeling impreciseness using IT2 FNs.
Therefore the present study can be extended by representing
impreciseness using general T2 FNs. Also, we have used conventional
multi-objective optimization techniques to deal with conflicting
objectives. So it is also a matter of further investigation to deal
with multiple objectives of the problem using evolutionary
algorithms like Multi-Objective Genetic Algorithm (MOGA) and
Non-dominated Sorting Genetic Algorithm (NSGA).

\subsection*{Acknowledgements:} The authors are thankful to the Editor and the anonymous Reviewers for valuable suggestions which lead to an improved version of the manuscript.

\end{document}